\documentclass[12pt]{amsart}
\usepackage{amsfonts,amssymb,amsthm,eucal}
\usepackage{xy}
\input xy
\xyoption{all}

\newtheorem{thm}{Theorem}
\newtheorem{cor}[thm]{Corollary}
\newtheorem{lemma}[thm]{Lemma}
\newtheorem{prop}[thm]{Proposition}

\newcommand{\inprod}[2]{\left\langle #1, #2 \right\rangle}

\DeclareMathOperator{\Cov}{Cov}
\DeclareMathOperator{\Var}{Var}
\DeclareMathOperator{\vol}{vol}
\DeclareMathOperator{\conv}{conv}
\renewcommand{\Re}{\operatorname{Re}}
\renewcommand{\Im}{\operatorname{Im}}
\DeclareMathOperator{\Id}{Id}

\newcommand{\Prob}{\mathbb{P}}
\newcommand{\Q}{\mathbb{Q}}
\newcommand{\R}{\mathbb{R}}
\newcommand{\E}{\mathbb{E}}

\newcommand{\Law}{\mathcal{L}}
\newcommand{\C}{\mathbb{C}}
\newcommand{\X}{\mathcal{X}}
\newcommand{\F}{\mathcal{F}}
\renewcommand{\epsilon}{\varepsilon}

\newcommand{\nought}{{\rm o}}

\newcommand{\ds}{\displaystyle}

\allowdisplaybreaks[3]

\title{The central limit problem for random vectors with symmetries}

\author[E.\ Meckes]{Elizabeth S.\ Meckes}
\email{meckes@math.stanford.edu}

\author[M.\ Meckes]{Mark W.\ Meckes}
\email{mark@math.stanford.edu}

\address{Department of Mathematics, Stanford University, Stanford,
California 94305, U.S.A.}

\begin{document}

\begin{abstract}
Motivated by the central limit problem for convex bodies, we study normal 
approximation
of linear functionals of high-dimensional random vectors with various
types of symmetries. In particular, we obtain results for distributions which
are coordinatewise symmetric, uniform in a regular simplex, or spherically
symmetric.  Our proofs are based on Stein's method of exchangeable pairs; 
as far as we know, this approach has not previously been used in convex
geometry and we give a brief introduction to the classical method.
The spherically symmetric case is treated by a variation 
of Stein's method which is adapted for continuous symmetries.
\end{abstract}

\maketitle

\section{Introduction}

Given a random vector $X = (X_1, \ldots, X_n) \in \R^n$, $n\ge 2$, and a
fixed $\theta \in S^{n-1}$, consider the random variable
$$W_\theta = \inprod{X}{\theta},$$
where $\inprod{\cdot}{\cdot}$ is the standard inner product on $\R^n$.
A typical example of interest is when $X$ is distributed uniformly
in a convex body.
In this paper we are interested in
determining sufficient conditions under which 
$W_\theta$ is approximately
normal, and in obtaining specific error estimates, possibly
depending on $\theta$.  To do this, we apply Stein's method of 
exchangeable pairs.  This technique has not previously been used in studying
problems from convex geometry, and we believe it will continue to be 
useful in that context.

To begin with, we will assume that $X$ is {\em isotropic}, that is,
that $W_\theta$ has mean $0$ and variance $1$ for every $\theta \in S^{n-1}$.
Equivalently, $X$ is isotropic if 
$$ \E X_i = 0 \quad \mbox{and} \quad \E X_i X_j = \delta_{ij}, $$
where $\delta_{ij}$ is the Kronecker delta. This is no real restriction,
since every random vector with finite second moment which is not supported
on a proper affine subspace has an affine image which is isotropic.

In the case that the components of $X$ are independent, bounds on the
distance of $W_\theta$ from normal follow from classical results.
For example, the Berry-Esseen theorem for sums of independent, non-identically
distributed random variables implies that if $X$ is isotropic with 
independent components, then
$$\sup_{t\in\R}\big| \Prob[W_\theta \le t] - \Phi(t) \big|
  \le 0.8 \Big(\max_{1\le i \le n} \E |X_i|^3 \Big) 
  \sum_{j=1}^n |\theta_j|^3,$$
where $\Phi$ denotes the standard normal distribution function. There is 
also a body of work going back to Sudakov \cite{Sudakov} 
(see \cite{Bobkov1} for a recent contribution and further references) 
on randomized central limit theorems, which show that under quite 
general assumptions
on $X$,  $W_\theta$ is close to some average distribution for most 
$\theta$.  In these results the average distribution
may or may not be normal, and ``most'' may refer to the rotation invariant
probability measure on $S^{n-1}$ or to some other distribution on weights.

Our motivation in studying this problem comes in part from the so-called
central limit problem for convex bodies, which is to show that the uniform
measure on any 
high-dimension\-al convex body has some one-dimensional projection which is
approximately normal, or that most one-dimensional projections 
are approximately normal.  Most of the results on this problem 
\cite{ABP,BHVV,Sodin}
prove some form of the latter conjecture (under appropriate assumptions
on the convex body), and thus fit into the framework of randomized central
limit theorems; none of the results in these
papers identify any specific $\theta$ for which $W_\theta$ is approximately
normal. The paper \cite{BV}
studies approximate normality of $W_\theta$ for specific $\theta$ 
when $X$ is uniformly distributed
in a cube, Euclidean ball, crosspolytope, or simplex, but in the last two
cases only for a very restricted set of $\theta$ and with rather limited
quantitative information.

Of course, there is no hope to identify any specific $\theta$ for which
$W_\theta$ is approximately normal without some additional assumptions on 
the distribution of $X$. Here the additional assumptions we consider are
more geometric than probabilistic in nature. Specifically, we consider
distributions which have a sufficiently rich class of symmetries, although
we emphasize that our results do not require $X$ to be drawn from a convex
body, or even to be continuous.  Stein's method, described in 
Section \ref{S:background}, allows us to take advantage of these symmetries
in order to reduce normal approximation to estimation of certain low-order 
moments.

Our first main result 
treats distributions which are symmetric with respect to reflection in a
suitable collection of hyperplanes. 
Our hypothesis encompasses both the class of distributions which are
coordinatewise symmetric (Corollary \ref{T:uncon}) and those with the 
symmetries
of a regular simplex (Corollary \ref{T:simplex}). The error
bounds for the approximations are in many cases small enough to derive
multivariate randomized versions which improve on
existing results.

Our second main result, treating spherically symmetric distributions,
is proved by
a variation of the classical version of Stein's method adapted to
take advantage of continuous symmetries.  This result has as corollaries
several classical results as well as some new applications.  We also
make connections with Poincar\'{e} inequalities for probability measures on
$\R^n$.

The layout of this paper is as follows. We first define notations which
will be used throughout the paper. In Section \ref{S:statements} we
state our results and several corollaries, and give comparisons to existing 
results. Section \ref{S:background} gives a brief introduction to 
Stein's method.  
Section \ref{S:proofs-1} contains the proofs of the first main
theorem 
and its corollaries. Section \ref{S:rotations} contains the proof of
the second main theorem using the variation of Stein's method described
above, and the proofs of its corollaries. 

\subsection*{Notation}

Let $\ell_p^n = (\R^n, \| \cdot \|_p)$, where $\| \cdot \|_p$ denotes the norm
$$\| x \|_p = \left( \sum_{i=1}^n |x_i|^p \right)^{1/p}$$
for $1\le p < \infty$, and 
$$\| x \|_\infty = \max_{1\le i \le n} |x_i|.$$

For $v \in R^n$, define $v \otimes v:\R^n \to \R^n$ by
$$v\otimes v(x)=\inprod{x}{v}v;$$
if $v\in S^{n-1}$, $v\otimes v$ is the orthogonal projection onto
the span of $v$.
A set of vectors $u_1, \ldots, u_m \in S^{n-1}$ such that
\begin{equation}\label{E:utf}
\sum_{i=1}^m u_i \otimes u_i = \frac{m}{n} I_n,
\end{equation}
is known in the signal processing literature as a {\em normalized tight frame} 
($I_n$ is the identity on $\R^n$). By taking the trace of both sides of
\eqref{E:utf}, one can see that $\frac{m}{n}$ is the only possible constant
that can appear.

The Grassmann manifold of $k$-dimensional subspaces of $\R^n$ is denoted
$G_{n,k}$; it is equipped with a unique rotation-invariant probability measure
$\lambda_{n,k}$.
For a fixed subspace $E \subset \R^n$, let $P_E$ denote the orthogonal 
projection
onto $E$, and $\gamma_E$ the standard Gaussian measure on $E$. 

For a random variable or random vector $X$, let $\Law(X)$ denote the
distribution of $X$
Given two probability measures $\mu$ and $\nu$ on $E$, define the
$T$-distance between them as
$$T(\mu,\nu) = \sup \bigl\{ |\mu(H) - \nu(H)| : H \mbox{ is an affine 
  half-space of $E$} \bigr\}. $$
In particular,
$$T\bigl(\Law(P_E(X)), \gamma_E\bigr) 
  = \sup_{\theta \in E \cap S^{n-1}, \ t\in \R}
  \bigl| \Prob[W_\theta \le t] - \Phi(t) \bigr|.$$
This is a geometrically natural notion of distance between measures since it
is invariant under nonsingular affine transformations and is thus not tied
to any particular coordinate system.  In addition, the topology induced by
$T$ on the space of probability measures on $E$ is stronger than the
$w^*$ topology.

The total variation distance between two probability measures $\mu$ and $\nu$
is
$$d_{TV}(\mu, \nu) = 2 \sup \bigl\{|\mu(A) - \nu(A)| : A \mbox{ is measurable}
  \bigr\}.$$
Recall that if $\mu$ and $\nu$ both have densities, then their total variation
distance is the $L_1$ distance between their densities.

Finally, note that symbols like $c$, $c_1$, etc.\ which represent absolute
constants may have different values from one appearance to the next.

\section{Statements of results} \label{S:statements}

Theorem \ref{T:proj}, the first main result of this paper, is based on 
existing normal approximation results proved via Stein's method. Corollaries
\ref{T:uncon}, \ref{T:SNCP}, \ref{T:lp}, and \ref{T:simplex} 
are all applications of Theorem \ref{T:proj}.

\begin{thm}\label{T:proj}
Let $u_1,\ldots,u_m \in S^{n-1}$ be a normalized tight frame, and for any
$x\in \R^n$ let $x_{(i)} = \inprod{x}{u_i}$.
Suppose that $X$ is a random vector whose distribution is
invariant under reflections in each of the hyperplanes $u_i^\perp$.  
Let $\theta\in S^{n-1}$ be fixed. Then
\begin{equation*}
\begin{split}
\sup_{t\in\R} \bigl|\Prob[W_\theta\le t]-\Phi(t)\bigr|
& \le 2\sqrt{\frac{n^2}{m^2} \sum_{i,j=1}^m \theta_{(i)}^2 \theta_{(j)}^2
  \E \bigl[X_{(i)}^2 X_{(j)}^2 \bigr] - 1} \\
&\quad + \left(\frac{8}{\pi}\right)^{1/4} \sqrt{\frac{n}{m}
  \Big(\max_{1\le j \le m} \E |X_{(j)}|^3\Big) 
  \sum_{i=1}^m |\theta_{(i)}|^3}.
\end{split}
\end{equation*}
If in addition $\ds\max_{1 \le i \le m} |X_{(i)}| \le a$ almost surely, then
\begin{equation*}
\begin{split}
\sup_{t\in\R} \bigl|\Prob[W_\theta \le t]-\Phi(t)\bigr|
& \le 24\sqrt{\frac{n^2}{m^2} \sum_{i,j=1}^m \theta_{(i)}^2 \theta_{(j)}^2
  \E \bigl[X_{(i)}^2 X_{(j)}^2 \bigr] - 1} \\
& \quad + 172na^3\max_{1 \le i \le m}|\theta_{(i)}|^3.
\end{split}
\end{equation*}
\end{thm}

\bigskip

The constants that appear in these statements, and explicit constants which
appear in any of the results that follow, are generally not the best 
possible, and are included only for concreteness.

It is not obvious from the statement that the error estimate in Theorem
\ref{T:proj} is useful.  However, the proofs of Corollaries \ref{T:uncon}
and \ref{T:simplex} will show that Theorem \ref{T:proj} allows different
cases of geometric interest to be treated easily in this unified
framework.

\medskip

Borrowing terminology from the geometry of Banach spaces, call
$X$ {\em unconditional} if its distribution is invariant under reflections 
in the coordinate hyperplanes, or equivalently, if $X$ has the same 
distribution as $(\epsilon_1 X_1, \ldots, \epsilon_n X_n)$ for any choice of
$\epsilon \in \{-1, 1\}^n$.

\begin{cor}\label{T:uncon}
Let $X$ be unconditional and isotropic, and let $\theta\in S^{n-1}$ be 
fixed. Then
\begin{equation*}
\begin{split}
\sup_{t\in\R} \bigl|\Prob[W_\theta \le t] - \Phi(t) \bigr| \le & \ 
  2 \sqrt{\Big(\max_{1\le i \le n} \E X_i^4 \Big)
    \| \theta \|_4^4 + \max_{i\neq j} \Cov(X_i^2, X_j^2)} \\
& + \left(\frac{8}{\pi}\right)^{1/4}
    \left(\max_{1\le i \le n}\sqrt{\E|X_i|^3}\right) \| \theta \|_3^{3/2}.
\end{split}
\end{equation*}
If moreover $X \in [-a,a]^n$ almost surely, then
\begin{equation*}
\begin{split}
\sup_{t\in\R} \bigl|\Prob[W_\theta \le t] - \Phi(t) \bigr| \le & \ 
  24 \sqrt{\Big(\max_{1\le i \le n} \E X_i^4 \Big) 
    \| \theta \|_4^4 + \max_{i\neq j} \Cov(X_i^2, X_j^2)} \\
& + 172 n a^3 \| \theta\|_\infty^3.
\end{split}
\end{equation*}
\end{cor}

\bigskip

The statement of Corollary \ref{T:uncon} motivates the introduction of
the following definition, taken from \cite{NR}. A random vector $X$ has the
{\em square negative correlation property} if 
$$\E X_i^2 X_j^2 \le (\E X_i^2) (\E X_j^2) \quad \mbox{for $i\neq j$},$$
i.e.\ if $\Cov(X_i^2, X_j^2) \le 0$ for $i \neq j$.

A random vector $X$ is called {\em log-concave} if it has a logarithmically
concave density with respect to Lebesgue measure on $\R^n$.  In particular,
a random vector which is uniformly distributed on a convex body is
log-concave.

It natural to conjecture (cf.\ Section 5 of \cite{NR}) that an 
isotropic unconditional
log-concave random vector has the square negative correlation property.
However, this is not the case, as shown by the example \cite{Bobkov3} of
the density
\begin{equation} \label{E:Bobkov}
a_n e^{-b_n \| x \|_\infty}
\end{equation}
on $\R^n$, where $a_n$ and $b_n$ are appropriate normalizing constants.
(Counterexamples also exist which are uniformly distributed in a convex
body \cite{Bobkov3}.) However, the weaker conjecture that under these 
conditions 
\begin{equation} \label{E:weak-correlation}
\Cov (X_i^2, X_j^2) \le \frac{c}{n} \quad \mbox{for $i\neq j$},
\end{equation}
for some absolute constant $c$, is open (and is in particular satisfied
by $X$ with the density \eqref{E:Bobkov}). This conjecture is related 
to the Kannan-Lov\'asz-Simonovits conjecture on isoperimetric constants
\cite{KLS} (see \cite{BK} for a discussion of this issue, and cf.\
Corollary \ref{T:Poincare} below).

The error bounds in Theorem \ref{T:proj} are small enough in many cases
to show that $W_\theta$ is uniformly close to normal for all unit vectors
$\theta$ in a typical subspace $E\subset \R^n$ of relatively large dimension. 
This means that the projection of $X$ onto $E$ is close to normally
distributed in the sense of $T$-distance. In order
to quantify this phenomenon, for $\epsilon>0$, define
$$A_{n,k}(\epsilon) = \bigl\{ E\in G_{n,k} : 
  T\bigl(\Law(P_E(X)), \gamma_E\bigr) \le \epsilon \bigr\};$$
that is, $A_{n,k}(\epsilon)$ is the set of $k$-dimensional subspaces 
$E \subset \R^n$ such that the projection of $X$ onto $E$ is $\epsilon$-close 
to normal in the sense of $T$-distance.

The following lemma, proved in \cite{BHVV}, allows normal approximation
in the sense of total variation to be deduced from approximation of 
distribution functions in the case of log-concave distributions.

\begin{lemma}[Brehm-Hinow-Vogt-Voigt] \label{T:tv}
There is an increasing function $\beta:(0,\infty) \to (0,2]$ satisfying
$$\beta(t) = O\Big(\sqrt{-t\log t}\Big) \quad \mbox{as $t\to 0$}$$
such that for any log-concave random variable $W$,
$$d_{TV}\big(\Law(W),\gamma_\R\big) \le 
  \beta\left(\sup_{t\in\R}\big| \Prob[W\le t] - \Phi(t) \big|\right).$$
\end{lemma}

\bigskip

All references to $\beta$ in the next three results are to the function
in Lemma \ref{T:tv}.

\begin{cor}\label{T:SNCP}
Let $X$ be unconditional and isotropic, with the square negative 
correlation property. Then there are constants $c_1, \ldots, c_{11}$,
independent of $X$ and $n$, such that each of the following holds.
\begin{enumerate}
\item \label{I:moment-bounds} 
If $\E |X_i|^3 \le a$ and $\E X_i^4 \le b$ for all $i$, then
$$\sup_{t\in\R} \bigl|\Prob[W_\theta \le t] - \Phi(t) \bigr|
  \le 2\big(\sqrt{b}\| \theta \|_4^2 +
  \sqrt{a} \| \theta \|_3^{3/2}\bigr)$$
for all $\theta\in S^{n-1}$, and
$$\lambda_{n,k} \big(A_{n,k}(\epsilon)\big)
  \ge 1 - c_1 \exp\left[-c_2 \min\{a^{-2/3}\epsilon^{4/3}, b^{-1/2}\epsilon\} n
  \right]$$
for $\epsilon \ge c_3 \max\{\sqrt{a}n^{-1/4}, \sqrt{b} n^{-1/2}\}$ and 
$k \le c_4 \min\{a^{-2/3}\epsilon^{4/3}, b^{-1/2}\epsilon\} n$.

\item \label{I:log-concave}
If $X$ is log-concave, then
$$\sup_{t\in\R} \bigl|\Prob[W_\theta \le t] - \Phi(t) \bigr|
  \le c_5 \| \theta \|_3^{3/2}$$
for all $\theta\in S^{n-1}$, and
$$\lambda_{n,k} \big(A_{n,k}(\epsilon)\big)
  \ge 1 - c_1 \exp\left[-c_6 \epsilon^{4/3} n\right]$$
for $\epsilon \ge c_7 n^{-1/4} $ and $k \le c_8 \epsilon^{4/3} n$.
Furthermore,
$$d_{TV}\big(\Law(W_\theta),\gamma_\R\big) \le 
  \beta\big(c_5 \| \theta \|_3^{3/2}\big)$$
for all $\theta \in S^{n-1}$.

\item \label{I:bounded}
If $X \in [-a,a]^n$ almost surely, then
$$\sup_{t\in\R} \bigl|\Prob[W_\theta \le t] - \Phi(t) \bigr|
  \le 196 n a^3 \| \theta \|_\infty^3$$
for all $\theta\in S^{n-1}$, and
$$\lambda_{n,k} \big(A_{n,k}(\epsilon)\big)
  \ge 1 - c_1 \exp\left[-c_{9} a^{-2} \epsilon^{2/3} n^{1/3}\right]$$
for $\epsilon \ge c_{10} a^3 (\log n)^{3/2} n^{-1/2}$ and 
$k \le c_{11} a^{-2} \epsilon^{2/3} n^{1/3}$.
If moreover $X$ is log-concave, then
$$d_{TV}\big(\Law(W_\theta),\gamma_\R\big) \le 
  \beta\big(196 n a^3 \| \theta \|_\infty^3\big)$$
for all $\theta \in S^{n-1}$.
\end{enumerate}
\end{cor}

\bigskip

The square negative correlation property is included as a hypothesis
of Corollary \ref{T:SNCP} only for convenience.  Replacing it with the
hypothesis \eqref{E:weak-correlation} would result only in a
weakening of the constants that appear, and in fact the even weaker
hypothesis
$$\Cov (X_i^2, X_j^2) \le \frac{c}{\sqrt{n}} \quad \mbox{for $i\neq j$}$$
would suffice for the same conclusion in part \ref{I:log-concave}.
In particular, the conclusion of part \ref{I:log-concave} applies to
$X$ distributed according to the density \eqref{E:Bobkov}.

One could also deduce randomized total variation results for one-dimensional
projections in parts \ref{I:log-concave} and \ref{I:bounded} of
Corollary \ref{T:SNCP}, but the statements are more complicated.

Naor and Romik \cite[Theorem 5]{NR} proved a result comparable to the
randomized statement in part \ref{I:moment-bounds} of Corollary 
\ref{T:SNCP}.  Under similar hypotheses (but without unconditionality),
they showed
$$\lambda_{n,k} \big(A_{n,k}(\epsilon)\big)
  \ge 1 - \frac{c_1}{\epsilon} \exp\left[-c_2 b^{-1} \epsilon^4 n\right]$$
for $\epsilon > 0$ and $k \le c_3 b^{-1} \epsilon^4 n$, so Corollary
\ref{T:SNCP} improves on the dependence on both $b$ and $\epsilon$ in
the unconditional case. 
In the case that $X$ is uniform in a convex body (hence log-concave)
and has the square negative correlation property, Antilla, Ball, and
Perissinaki \cite[Theorem 2]{ABP} showed
\begin{equation}\label{E:ABP}
\lambda_{n,1} \big(A_{n,1}(\epsilon) \big)
  \ge 1 - n \exp \left[-c_1 \epsilon^2 n\right]
\end{equation}
for $\epsilon \ge c_2 n^{-1/3}$. Part \ref{I:log-concave} of 
Corollary \ref{T:SNCP} improves (in the unconditional case) 
on their dependence on $\epsilon$, although
for a slightly more restricted range of $\epsilon$,
and does not require that $X$ be chosen from a convex body.
In the case of certain bounded distributions,
part \ref{I:bounded} improves further on the results of \cite{ABP,NR},
as will be illustrated by Corollary \ref{T:lp} below.

\medskip

The next corollary treats a class of examples of particular interest 
in asymptotic convex geometry, namely, $X$ chosen
from various natural distributions on the unit balls of the spaces $\ell_p^n$.
In addition to the uniform measure on the interior, there are two 
geometrically natural measures on the boundary of a convex
body $K$ whose interior contains the origin. First there is
$(n-1)$-dimensional Hausdorff measure, or surface measure. Second, there
is cone measure $\mu$, defined by
$$\mu(A) = \vol \Big( \bigcup_{t\in [0,1]} tA \Big) \quad
  \mbox{for $A \subset \partial K$}.$$
Cone measure is the measure on $\partial K$ for which there is
a straightforward extension of the familiar polar integration formula,
with $\partial K$ replacing $S^{n-1}$.

\begin{cor}\label{T:lp}
For $1\le p\le \infty$, let $X$ have one of the following distributions:
\begin{enumerate}
\item uniform measure on the ball of $\ell_p^n$, scaled to be isotropic;
\item \label{I:cone}
normalized cone measure on the sphere of $\ell_p^n$, scaled to be 
isotropic; or
\item \label{I:surface}
normalized surface measure on the sphere of $\ell_p^n$, scaled such
that the normalized cone measure is isotropic.
\end{enumerate}
Then there are absolute constants $c_1, \ldots, c_5$ and constants 
$d_{1,p},\ldots,d_{4,p}$ depending only on $p$ such that
\begin{equation}\label{E:lp}
\sup_{t\in\R} \bigl|\Prob[W_\theta \le t] - \Phi(t) \bigr|
  \le \min \left\{ c_1 \| \theta \|_3^{3/2}, 
       d_{1,p} n^{1+\frac{3}{p}} \| \theta \|_\infty^3\right\}
\end{equation}
for all $\theta\in S^{n-1}$,
$$\lambda_{n,k} \big(A_{n,k}(\epsilon)\big)
  \ge 1 - c_2 \exp\left[-c_3 \epsilon^{4/3} n\right]$$
for $\epsilon \ge c_4 n^{-1/4}$ and $k \le c_5 \epsilon^{4/3} n$,
and
$$\lambda_{n,k} \big(A_{n,k}(\epsilon)\big)
  \ge 1 - c_2 \exp\left[-d_{2,p} \epsilon^{2/3} n^{\frac{1}{3}-\frac{2}{p}}
  \right]$$
for $\epsilon \ge d_{3,p} (\log n)^{3/2} n^{\frac{3}{p}-\frac{1}{2}}$ and 
$k \le d_{4,p} \epsilon^{2/3} n^{\frac{1}{3} - \frac{2}{p}}$.

Furthermore, in the case that $X$ is chosen uniformly from the rescaled 
$\ell_p^n$ ball,
$$d_{TV}\big(\Law(W_\theta),\gamma_\R\big) \le 
  \beta\left(\min \left\{ c_1 \| \theta \|_3^{3/2}, 
       d_{1,p} n^{1+\frac{3}{p}} \| \theta \|_\infty^3\right\}\right)$$
for all $\theta \in S^{n-1}$.
\end{cor}

\bigskip

To compare the two bounds in \eqref{E:lp}, note that since $\|\theta\|_\infty
\le \|\theta \|_3 \le \|\theta\|_2 = 1$ and $\|\theta\|_3 \ge n^{-1/6}$,
the first error bound is better for all $\theta$ when $1\le p < 4$
(ignoring the constant factors). On the
other hand, for the principal diagonal $\theta = n^{-1/2}\sum e_i$ (which
roughly captures typical behavior for $\ell_p^n$ norms on $S^{n-1}$) 
the second bound is better for $p \ge 12$. 
In particular, for $p>18$, Corollary \ref{T:lp} improves on the
typical rate of convergence to normality of about $n^{-1/3}$ which follows 
from \eqref{E:ABP}. Corollary \ref{T:lp} also improves on Theorems 7 and 8 of
\cite{NR}, which show
$$\lambda_{n,k} \big(A_{n,k}(\epsilon)\big)
  \ge 1 - \frac{c_1}{\epsilon} \exp\left[-c_2 \epsilon^4 n\right]$$
for $\epsilon > 0$ and $k \le c_3 \epsilon^4 n$ in cases \ref{I:cone}
and \ref{I:surface} of Corollary \ref{T:lp}.

Brehm and Voigt \cite{BV} considered $X$ uniformly distributed in the
rescaled $\ell_p^n$ ball for $p=1,2,\infty$.  See the discussion of
Corollary \ref{T:Euclidean} below for the case $p=2$.  
In the case $p=\infty$, they
derive sharper error bounds for general $\theta$ than here. In the case
$p=1$ (the crosspolytope) however, they consider only the case
$\theta = n^{-1/2} \sum e_i$ for $n \to \infty$, and do not obtain an 
explicit rate of convergence, so Corollary \ref{T:lp} provides 
a substantial generalization and strengthening.

\medskip

As discussed earlier, the form of Theorem \ref{T:proj} is 
general enough to accommodate the symmetries both of unconditional 
distributions
and of a regular simplex, which is treated in the next result.  
\begin{cor}\label{T:simplex}
Let $\Delta_n=\sqrt{n(n+2)}\conv\{v_1,\ldots,v_{n+1}\}$ be an isotropic 
regular simplex, where $v_i \in S^{n-1}$.  
Let $X$ be uniformly distributed in $\Delta_n$, and let
$\theta \in S^{n-1}$ be fixed. Then there are constants $c_1,\ldots, c_5$,
independent of $n$, such that
\begin{equation*}
\sup_{t\in\R} \bigl|\Prob[W_\theta \le t] - \Phi(t) \bigr|
\le c_1 \sqrt{\sum_{i=1}^{n+1}\big|\inprod{\theta}{v_i} \big|^3}
\end{equation*}
and
$$d_{TV}\big(\Law(W_\theta),\gamma_\R\big) \le 
  \beta\left(c_1 \sqrt{\sum_{i=1}^{n+1}\big|\inprod{\theta}{v_i} \big|^3}
  \right)$$
for all $\theta\in S^{n-1}$, and 
$$\lambda_{n,k} \big(A_{n,k}(\epsilon)\big)
  \ge 1 - c_2 \exp\left[-c_3 \epsilon^{4/3} n\right]$$
for $\epsilon \ge c_4 n^{-1/4}$ and $k \le c_5 \epsilon^{4/3} n$.
\end{cor}

\bigskip

The case in which $X$ is uniformly distributed in a regular simplex was
also considered in \cite{BV}. However, the results there consider only
a certain discrete set of $\theta$ (roughly those for which the behavior
of $W_\theta$ is best) for $n \to \infty$, and do not derive any explicit 
rate of convergence to normality.

\medskip

The remaining results are not based on existing normal approximation
theorems; instead, the proofs use a variation of Stein's method of 
exchangeable pairs, adapted to situations in which there are continuous
symmetries.  This variation was introduced
by Stein in \cite{Stein2} and developed further by the first-named author in
\cite{EMeckes} in studying functions on the classical matrix groups. 

\begin{thm}\label{T:sph-symm}
Let $X$ be an isotropic random vector, with finite third moment, whose 
distribution is spherically symmetric. Then for any $\theta \in S^{n-1}$,
\begin{equation*}\begin{split}
d_{TV}\bigl(\Law(W_\theta), \gamma_\R \bigr) 
  &\le 4 \E \bigl|1 - \E[X_2^2 | X_1]\bigr|\\
&\le \frac{4}{n-1}
\E\big|\|X\|_2^2 - n \big|+\frac{8}{n-1}\\
&\le \frac{4}{n-1}
\sqrt{\Var\left(\|X\|_2^2\right)}+\frac{8}{n-1}.\end{split}\end{equation*}
\end{thm}

\bigskip

The latter bounds in Theorem \ref{T:sph-symm} reduce normal approximation of
one-dimensional projections of spherically symmetric random vectors to
the problem of estimating deviations of $\|X\|_2^2$ from its mean. 
Some kind of such concentration of $\| X \|_2$ is either explicitly a 
hypothesis, or closely related to the
key hypothesis, in many of the existing results on the central limit problem
for convex bodies, cf.\ \cite{ABP,BK,BHVV,Sodin}.

Before proceeding to some classical consequences of Theorem 
\ref{T:sph-symm}, we first state a corollary giving a new
randomized central limit theorem, which in particular gives information about
the central limit problem for convex bodies.
Suppose $X$ is a (not necessarily spherically symmetric) isotropic random
vector and $U$ is a random $n\times n$ orthogonal matrix, distributed
according to Haar measure independently of $X$. 
Notice that the distribution of $W = \inprod{UX}{e_1} = \inprod{X}{U^{-1}e_1}$ 
is the average
(with respect to the rotation invariant probability on $S^{n-1}$) of the
distributions of $W_\theta$ over all $\theta\in S^{n-1}$. As mentioned in
the introduction, the distribution of $W$ is the object of some of the
work on randomized central limit theorems. Since $\tilde{X} = UX$ is
a spherically symmetric isotropic random vector and $\|\tilde{X}\|_2
= \|X\|_2$, the following is an immediate consequence of Theorem
\ref{T:sph-symm}.

\begin{cor}\label{T:avg-dist}
Let $X$ be an isotropic random vector with finite third moment, and let
$W$ be as defined above. Then
\[\begin{split}
d_{TV} \big( \Law(W), \gamma_\R \big) & \le
\frac{4}{n-1} \E \big| \|X\|_2^2 - n \big| + \frac{8}{n-1} \\
& \le \frac{4}{n-1} \sqrt{\Var\left(\|X\|_2^2\right)}+\frac{8}{n-1}.
\end{split}\]
\end{cor}

\bigskip

In the case in
which $X$ is uniformly distributed in $K \subset \R^n$, the density of
$W$ gives the average $(n-1)$-dimensional volume of a hyperplane section
of $K$ at a given distance from the origin. 
For $K$ a convex body, Bobkov and Koldobsky \cite{BK} proved a pointwise 
bound on the difference between the density of $W$ and the standard 
normal density, with a bound which also explicitly involves the variance
which appears in Corollary \ref{T:avg-dist} (and which is also essentially
of the order $n^{-1}$ as long as the variance is not too big). 
Corollary \ref{T:avg-dist} gives instead an $L_1$ bound on the difference of 
these densities. See \cite{KL} for an earlier asymptotic result in the case
that $K$ is a cube, and \cite{BB} for a generalization of the result
of \cite{BK} to arbitrary distributions and a multivariate version for 
sections by $k$-codimensional affine subspaces.

\medskip

The following easy corollary of Theorem \ref{T:sph-symm} is well-known; 
versions for higher-dimensional projections are proved in \cite{DF2} for the
sphere and in \cite{BV} for the ball.
Theorem \ref{T:sph-symm} allows the cases of both the ball and the
sphere to be presented simply as part of a unified framework.

\begin{cor}\label{T:Euclidean}
If $X$ has the uniform distribution on the Euclidean ball of radius
$\sqrt{n+2}$ or the uniform distribution on the sphere of radius $\sqrt{n}$,
then for any $\theta \in S^{n-1}$
$$d_{TV}\bigl(\Law(W_\theta), \gamma_\R \bigr) \le \frac{a}{n-1},$$
where $a=16$ in the case of the ball and $a=8$ in the case of the ball.
\end{cor}

\bigskip

Corollary \ref{T:Euclidean} gives the correct order of approximation in
both cases, although the constants can be improved.

The first error estimate in Theorem \ref{T:sph-symm} is also strong enough
to recover, as an immediate consequence,
a version of the characterization of the normal distribution
as the unique spherically symmetric product measure on $\R^n$.
(The two-dimensional version of this characterization is the classical 
Herschel-Maxwell theorem; see \cite{Bryc} for various other versions.)

\begin{cor}\label{T:Normal}
A random vector $X$ with finite third moment has the standard normal 
distribution if and only if
$X$ is spherically symmetric and has independent components with variance
$1$.
\end{cor}

\bigskip

Recall that a random vector $X$ is said to satisfy a {\em Poincar\'e
inequality} with constant $\lambda_1$ (the {\em spectral gap} of $X$) if
\begin{equation}\label{E:Poincare}
\lambda_1 \Var \big(f(X)\big) \le \E \|\nabla f(X)\|_2^2
\end{equation}
for every smooth $f:\R^n\to \R$. The last estimate in Theorem 
\ref{T:sph-symm} provides a connection between normal approximation and
spectral gap estimates.  A similar connection has been observed in 
a different but related context by Bobkov and Koldobsky \cite{BK}.

\begin{cor}\label{T:Poincare}
Let $X$ be an isotropic spherically symmetric random vector with spectral
gap $\lambda_1$. Then for any $\theta\in S^{n-1}$,
$$d_{TV}\big(\Law(W_\theta),\gamma_\R\big)
  \le \frac{10}{\sqrt{n \lambda_1}}.$$
\end{cor}

\bigskip

One concrete application of Corollary \ref{T:Poincare} is the following.
\begin{cor}\label{T:exp}
 Let 
$X$ have the isotropic spherically symmetric exponential density
$$a_n e^{-b_n \|x\|_2},$$
where $a_n$ and $b_n$ are appropriate normalization constants.
Then
$$d_{TV}\big(\Law(W_\theta),\gamma_\R\big) \le \frac{10\sqrt{13}}{n^{1/2}}.$$ 
\end{cor}

\bigskip

Bobkov \cite{Bobkov4} showed that for this distribution, 
$\frac{1}{13} \le \lambda_1 \le 1$, and so Corollary \ref{T:exp}
is immediate from Corollary \ref{T:Poincare}.  Since this 
distribution is given explicitly, one can obtain an 
error estimate of the same order by directly estimating
the variance in Theorem \ref{T:sph-symm};  however,
there is a large literature on spectral gap estimates in much 
less explicit contexts using only certain
geometric assumptions, typically diameter and/or curvature bounds
 (see \cite{Ledoux2} for a survey and further 
references).  Corollary \ref{T:Poincare} thus allows the treatment of 
distributions about which one has geometric information resulting in 
spectral gap estimates, but for which direct computation of the variance 
term in Theorem \ref{T:sph-symm} is not possible.

\medskip

There is also the following complex analogue of Theorem \ref{T:sph-symm}.
All of the previously defined notation of this paper used here should be
reinterpreted for vectors in $\C^n$ in the most obvious way.

\begin{thm}\label{T:unitary}
Let $X\in\C^n$ be a random vector with finite third moment 
such that $\E X_i=0$ for each $i$, 
$\E X_iX_j=\E X_i\overline{X}_j=0$ if $i\neq j$, and $\E(\Re X_i)^2=\E
(\Im X_i)^2=1$.  Suppose the distribution of $X$ is invariant under
multiplication by a unitary matrix.  Then for any $\theta\in S_{\C}^{n-1}$,
\begin{equation*}\begin{split}
d_{TV}  \bigl(\Law(\Re W_\theta), \gamma_\R \bigr) 
  &\le 4 \E \left|1 - \frac{n}{2(n+1)}\E\big[|X_2|^2 \big| X_1\big]\right|
+\frac{1}{n-1}\\
&\le\frac{2n}{n^2-1}\sqrt{\Var\left(\|X\|_2^2\right)}+\frac{5}{n-1}.
\end{split}\end{equation*}
\end{thm}

Note that under the hypotheses of Theorem \ref{T:unitary}, 
$\E [|X_2|^2 | X_1] = 2 \E[(\Re X_2)^2 | X_1]$, and so the appearance
of the factor of 2 inside the first bound above is to be expected.

\section{Background on Stein's method}
\label{S:background}

The essential idea of Stein's method is the notion of a characterizing 
operator.  Say that $T_\nought$ is a characterizing operator for a 
distribution $\mu$ on $\R$ if the following conditions hold:
\begin{enumerate}
\item $\int T_\nought f(t)\ d\mu(t)=0$ for all $f$ such that $T_\nought f$ is 
$\mu$-integrable, and 
\item if $\nu$ is a probability measure on $\R$ such that $\int T_\nought 
f(t) \ d\nu(t)=0$
for all $f$ with $T_\nought f$ $\nu$-integrable, then $\mu=\nu$.
\end{enumerate}
A characterizing operator is a strong
characterization of a distribution, in the sense that if $T_\nought$ is 
characterizing for $\mu$ and $\nu$ is a measure such that $\int T_\nought f(t)
\ d\nu(t)$
is small for a large class of test functions $f$, then $\nu\approx\mu$ 
in some sense. 

This idea is quantified in the method of exchangeable pairs as follows.
Let $W=W(\omega)$ be a random variable defined on a probability space 
$(\Omega, \Prob)$, 
and let $\X$ be a space of measurable functions on $\Omega$.  
Think of $\E$ as a linear map from $\X$ to $\R$, $\E f = \int f(\omega) \ 
d\Prob(\omega)$.  Let $\X_\nought$ be
a space of measurable functions on $\R$ and let $\E_\nought$ be the linear
function on $\X_\nought$ defined by $\E_\nought f=\int f(t) \ d\mu(t)$ for 
some fixed 
measure $\mu$.  In our applications, $\mu$ will be the standard normal
distribution, but one of the advantages of Stein's method is that 
the set-up is quite general and can be adapted to various other measures.
The random variable $W$ induces a map $\beta:\X_\nought\to\X$ defined by
\begin{equation}\label{beta}
\beta f(\omega)=f(W(\omega)).
\end{equation} 
Now, construct a symmetric probability $\Q$ on $\Omega\times\Omega$ with 
margins $\Prob$ (i.e., $\Q(A\times B)=\Q(B\times A)$ and $\Q(A\times\Omega)=
\Prob(A)$).  Note that this is the same as constructing an exchangeable pair
$(W,W') = (W(\omega), W(\omega'))$ from $W$.  
Let $\F$ be a space of measurable, antisymmetric
functions on $\Omega\times\Omega$ and use $\Q$ to define a map $T:\F\to\X$
by
\begin{equation}\label{T}
T f(\omega)=\E_\Q\left[f(\omega,\omega')\big|\omega\right];
\end{equation}
note that by exchangeability and antisymmetry, $\E T\equiv0$ on $\F$.
Let $\F_\nought$ be another space of measurable functions on $\R$, possibly the
same as $\X_\nought$, and let $T_\nought:\F_\nought\to\X_\nought$ be a 
characterizing operator of
$\mu$.  Let $U_\nought:\X_\nought\to\F_\nought$ be a pseudo-inverse to 
$T_\nought$, in the sense that 
\begin{equation}\label{T_o U_o}
T_\nought U_\nought f(t)=f(t)-\E_\nought f.  
\end{equation}
Finally, let $\alpha:\F_\nought\to\F$.  All of these definitions are summarized
in the following diagram:
\begin{equation} \label{diagram}
\xymatrix{\F\ar[rr]^{T}&&\X\ar[dr]^{\E}\\
&&&\R\\ \F_\nought\ar@<.5ex>[rr]^{T_\nought}\ar[uu]^\alpha&&\X_\nought
\ar@<.5ex>[ll]^{U_\nought}\ar[uu]^\beta\ar[ur]_{\E_\nought}}
\end{equation}

The following easy lemma \cite{Stein1}
is the quantitative version of the heuristic at the
beginning of the section.
\begin{lemma}[Stein]\label{L:stein}
Suppose that in the diagram of spaces and maps above, $\E T=0$ and 
$T_\nought U_\nought=\Id-\E_\nought$.  Then
\begin{equation*}
\E\beta-\E_\nought=\E(\beta T_\nought-T\alpha)U_\nought.
\end{equation*}
\end{lemma}

\bigskip
Note that for $f\in\X_\nought$, 
$$\E \beta f-\E_\nought f=\E f(W)-\E f(Z),$$
where $Z$ is a random variable with distribution $\mu$.  The strategy is 
apply Stein's lemma to bound this difference uniformly over a large
class of test functions.  

To apply the method of exchangeable pairs, one needs an approximating 
distribution and a characterizing operator.  In this paper, the 
measure $\mu$ will be standard Gaussian, with characterizing operator
$$T_\nought f(t)=f'(t)-tf(t)$$
and pseudo-inverse
$$U_\nought f(t)=e^{\frac{1}{2}t^2}
 \int_{-\infty}^t\big[f(s)-\E_\nought f\big]e^{-\frac{1}{2}s^2} \ ds.$$
That $\E_\nought T_\nought=0$ can be verified by integration by parts, 
and verifying 
$T_\nought U_\nought=\Id-\E_\nought$ is just calculus.

The following estimates for $U_\nought$ are proved in  
\cite[p.\ 25]{Stein1} and are useful in estimating the error term from
Lemma \ref{L:stein}:
\begin{align}
\label{E:U_o}
\| U_\nought f \|_\infty \le& \sqrt{\frac{\pi}{2}} \| f-\E_\nought f \|_\infty
  \le \sqrt{2 \pi} \| f\|_\infty, \\
\label{E:U_o'}
\| (U_\nought f)' \|_\infty \le& 2 \| f-\E_\nought f \|_\infty 
  \le 4 \| f \|_\infty, \\
\label{E:U_o''}
\| (U_\nought f)'' \|_\infty \le& 2 \| f' \|_\infty.
\end{align}

Using this general set-up one can prove the following abstract normal
approximation theorem, which is the main tool used in the proof of 
Theorem \ref{T:proj}.  The first statement
was proved by Stein in \cite[Lecture III]{Stein1}; the second statement 
was proved by Rinott and Rotar in \cite{RR}.

\begin{prop}[Stein; Rinott-Rotar]\label{T:Stein-RR}
Let $(W,W')$ be an exchangeable pair of random variables such
that
$$\E W = 0, \quad \E W^2 = 1,$$
and
$$\E[W-W'|W] = \lambda W$$
for some $\lambda \in (0,1)$. Then
$$
\sup_{t\in\R} \bigl| \Prob[W\le t] - \Phi(t) \bigr| \le
  \frac{1}{\lambda} \sqrt{\Var \E [(W-W')^2|W]}
  + (2\pi)^{-1/4} \sqrt{\frac{1}{\lambda}\E |W-W'|^3}.
$$
If moreover $|W-W'|$ is almost surely bounded, then
$$
\sup_{t\in\R} \bigl| \Prob[W\le t] - \Phi(t) \bigr| \le
  \frac{12}{\lambda} \sqrt{\Var \E [(W-W')^2|W]}
  + \frac{43}{\lambda}\|W-W'\|_\infty^3.
$$
\end{prop}

\bigskip

{\bf Remarks}
\begin{enumerate}
\item For Stein's lemma to hold, no assumptions are needed on the map
$\alpha:\F_\nought\to\F$ other than that it does in fact produce antisymmetric
functions.  For normal approximation, one usually uses the map
$$\alpha f(\omega, \omega')=a\big(W(\omega')-W(\omega)\big)\big[f(W(\omega'))+
f(W(\omega))\big],$$
for some suitable choice of $a$.
\item In the description of the method, we have
been quite cavalier about exactly which function spaces should be used.
Of course, the choice of the space of test functions $\X_\nought$ determines 
the type of convergence; in practice, one generally fixes $\X_\nought$
 first and chooses
the remaining spaces in some way which guarantees that all of the maps do
fit into the diagram as shown.  The diagram and Stein's lemma 
are described here mainly 
as motivation for the approach taken in Section \ref{S:rotations} and are
often not used directly but as a guide for how to proceed.
\item Proposition \ref{T:Stein-RR} has been applied in many different
contexts.  In particular, Holmes and Reinert \cite{HR} use it
in analyzing statistics similar to our $W_\theta$ while studying process
approximation for the bootstrap.
\item This paper is the first that we know of to apply the exchangeable pairs
approach to Stein's method in convex geometry.  However, Reitzner 
\cite{Reitzner} used a rather different approach to Stein's method (based 
on dependency graphs) in proving central limit theorems for random polytopes.
\item We have described here how to apply Stein's lemma to normal 
approximation; application to Poisson approximation has also been extensively
developed, see \cite{CDM}.

\end{enumerate}

\section{Proof of Theorem \ref{T:proj} and its consequences}
\label{S:proofs-1}

\begin{proof}[Proof of Theorem \ref{T:proj}]
Let $X$ be a random vector invariant under reflections in the hyperplanes
defined by a normalized tight frame $u_1,\ldots, u_m$.
Define an exchangeable pair $(W,W')$ as in Proposition \ref{T:Stein-RR} as
follows. Let $I$ be chosen
uniformly from $\{1,\ldots,m\}$, independently of $X$, and define
$$X' = X - 2X_{(I)}u_{I},$$
i.e., $X'$ is obtained from $X$ by reflection in the
hyperplane $u_I^\perp$.
Then $(X,X')$ is an exchangeable pair of random vectors by assumption. 
Define $W=W_\theta=\inprod{X}{\theta}$
and $W'=\inprod{X'}{\theta}$.
Now $\E W = 0$ and $\E W^2 = 1$ since $X$ is isotropic, and
\begin{equation*}
\begin{split}
\E [W-W' | W] &= 
  \E \left[\left.\frac{1}{m}\sum_{i=1}^m 2 X_{(i)}\theta_{(i)}
  \right| W \right] \\
&=\frac{2}{m}\E\left[\left.\inprod{\sum_{i=1}^m u_i\otimes u_i(X)}{\theta}
  \right|W\right] \\
&= \frac{2}{n}W
\end{split}
\end{equation*}
by equation \eqref{E:utf}.

To apply Proposition \ref{T:Stein-RR} (with $\lambda = \frac{2}{n}$),
it remains to estimate the quantities
$$\Var\E[(W-W')^2|W], \quad \E |W-W'|^3, \quad \mbox{and}
\quad \| W-W' \|_\infty$$
(the last in the case that $\max |X_{(i)}|\le a$ almost surely).
First,
$$\E \bigl(\E[(W-W')^2|W]\bigr) = \E \bigl(\E[W^2 + (W')^2 - 2WW'|W]\bigr)
  = \frac{4}{n},$$
and by the conditional form of Jensen's inequality,
\begin{equation*}
\begin{split}
\E \bigl(\E[(W-W')^2|W]\bigr)^2 &\le  \E \bigl(\E[(W-W')^2|X]\bigr)^2 \\
&= \E \bigl(\E [(2X_{(I)}\theta_{(I)})^2 | X ] \bigr)^2 \\
&=\frac{16}{m^2} \sum_{i,j=1}^m \theta_{(i)}^2 \theta_{(j)}^2 
  \E \bigl[X_{(i)}^2 X_{(j)}^2\bigr],
\end{split}
\end{equation*}
so
\begin{equation}\label{E:Var}
\begin{split}
\Var\E[(W-W')^2|W]
&= \E \bigl(\E[(W-W')^2|W]\bigr)^2 -\frac{16}{n^2} \\
&\le \frac{16}{n^2} \left(\frac{n^2}{m^2} \sum_{i,j=1}^m 
  \theta_{(i)}^2 \theta_{(j)}^2 \E \bigl[X_{(i)}^2 X_{(j)}^2 \bigr] - 1\right).
\end{split}
\end{equation}     

Next,
\begin{equation}\label{E:3rd-moment}
\begin{split}
\E |W-W'|^3 &= 8 \E | X_{(I)}\theta_{(I)}|^3 \\
&= \frac{8}{m} \sum_{i=1}^m |\theta_{(I)}|^3 \E |X_{(i)}|^3 \\
&\le \frac{8}{m} \left(\max_{1\le i\le m} \E |X_{(i)}|^3\right) 
  \sum_{i=1}^m|\theta_{(i)}|^3.
\end{split}
\end{equation}
Finally, if $\max | X_{(i)}| \le a$ almost surely, then
\begin{equation}\label{E:sup}
|W-W'| = 2|\theta_{(I)}||X_{(I)}| \le 2 a \max |\theta_{(i)}|
\end{equation}
almost surely.  
Inserting the estimates \eqref{E:Var}, \eqref{E:3rd-moment}, and
\eqref{E:sup} into Proposition \ref{T:Stein-RR} now
proves Theorem \ref{T:proj}.
\end{proof}

\bigskip
{\bf Remark:}
One can also prove, by the same method, 
a version of Theorem \ref{T:proj} for random vectors which are invariant
under reflections in subspaces of arbitrary dimensions. Rather than a
normalized tight frame, one would consider a decomposition
$$\sum_{i=1}^m P_{E_i} = \alpha I_n$$
for subspaces $E_1, \ldots, E_m \subset \R^n$. In applications this can
lead to a nontrivial approximation result if the maximum dimension
of the $E_i$ remains bounded (or grows very slowly) as $n$ grows. 
This context covers,
for example, $X$ uniformly distributed in the $\ell_p$-sum of $m$ copies
of a fixed symmetric convex body $K \subset \R^k$.

\medskip

\begin{proof}[Proof of Corollary \ref{T:uncon}]
The corollary follows easily from Theorem \ref{T:proj}, by taking as the
normalized tight frame the standard basis of $\R^n$ (so that $x_{(i)} = x_i$). 
The form of the second
error term in both statements follows immediately from Theorem \ref{T:proj}.
To estimate the first error term, note that
\begin{equation*}
\begin{split}
\sum_{i,j=1}^m \theta_i^2 \theta_j^2 \E \bigl[X_i^2 X_j^2\bigr]
&\le \left( \max_{1\le j \le m} \E X_j^4 \right) \sum_{i=1}^m \theta_i^4
  + \left( \max_{i\neq j} \E X_i^2 X_j^2 \right) 
  \sum_{i,j=1}^m \theta_i^2 \theta_j^2 \\
&= \left( \max_{1\le j \le m} \E X_j^4 \right) \sum_{i=1}^m \theta_i^4
  + \max_{i\neq j} \E X_i^2 X_j^2
\end{split}
\end{equation*}
since $\| \theta \|_2 = 1$, and 
$$\E X_i^2 X_j^2 = \Cov \bigl(X_i^2, X_j^2\bigr) + 1$$
since $X$ is isotropic.
\end{proof}

\bigskip

The following lemma is used to prove the randomized statements in 
Corollaries \ref{T:SNCP}, \ref{T:lp}, and
\ref{T:simplex}.  It follows easily from a concentration 
inequality implicitly proved by Gordon \cite{Gordon} 
(see Theorem 6 in \cite{NR} for an explicit statement),
together with the well-known asymptotic orders of the averages of
$\ell_p^n$ norms over $S^{n-1}$.

\begin{lemma}\label{T:Gordon}
There are absolute constants $c_1,\ldots,c_4$ such that the following
hold.
\begin{enumerate}
\item \label{I:Gordon-3}
If $\delta \ge c_1 n^{-1/6}$ and $k \le c_2 \delta^2 n$, then
$$\lambda_{n,k} (\{E\in G_{n,k} : \| \theta \|_3 \le \delta \ \forall
  \theta \in E \cap S^{n-1} \}) \ge 1 - c_3 e^{-c_4 \delta^2 n}.$$
\item \label{I:Gordon-4}
If $\delta \ge c_1 n^{-1/4}$ and $k \le c_2 \delta^2 n$, then
$$\lambda_{n,k} (\{E\in G_{n,k} : \| \theta \|_4 \le \delta \ \forall
  \theta \in E \cap S^{n-1} \}) \ge 1 - c_3 e^{-c_4 \delta^2 n}.$$
\item \label{I:Gordon-infty}
If $\delta \ge c_1 \sqrt{\frac{\log n}{n}}$ and $k \le c_2 \delta^2 n$, then
$$\lambda_{n,k} (\{E\in G_{n,k} : \| \theta \|_\infty \le \delta \ \forall
  \theta \in E \cap S^{n-1} \}) \ge 1 - c_3 e^{-c_4 \delta^2 n}.$$
\end{enumerate}
\end{lemma}

\bigskip

\begin{proof}[Proof of Corollary \ref{T:SNCP}] ~
\begin{enumerate}
\item
The first statement is immediate from the first statement of Corollary 
\ref{T:uncon}; the second follows from parts 
\ref{I:Gordon-3} and \ref{I:Gordon-4} of Lemma \ref{T:Gordon}.

\item It is a well-known consequence of Borell's lemma (cf.\ \cite[Section
2.2]{Ledoux1}) that there is an absolute constant $c$, independent of $X$, 
$n$, and $p$, such that
\begin{equation}\label{E:Borell}
\bigl(\E |\inprod{X}{y}|^p\bigr)^{1/p}
\le c p \bigl(\E |\inprod{X}{y}|^2\bigr)^{1/2}
\end{equation}
for any log-concave random vector $X$, $p\ge 2$, and fixed vector
$y$. Thus part \ref{I:moment-bounds}
applies with some absolute constants $a$ and $b$. Furthermore, since
$\|\theta\|_4 \le \|\theta\|_3 \le \|\theta\|_2 = 1$, the first term in the
r.h.s.\ of part \ref{I:moment-bounds} is not of larger order than the second
term.  The total variation bound follows from Lemma \ref{T:tv} and
the fact that any projection of a log-concave measure is again
log-concave \cite{Prekopa}.

\item From the second statement of Corollary \ref{T:uncon} and the
trivial estimate $\E X_i^4 \le a^4$ we obtain
$$\sup_{t\in\R} \bigl|\Prob[W_\theta \le t] - \Phi(t) \bigr|
  \le 24 a^2 \|\theta\|_4^2 + 172 n a^3 \| \theta \|_\infty^3.$$
The estimates
$$\| \theta \|_4 \le n^{1/4} \| \theta \|_\infty$$
and
$$1 = \|\theta\|_2 \le \sqrt{n} \|\theta\|_\infty$$
are well-known, and $a\ge 1$ since $X$ is isotropic. 
Therefore
$$a^2 \|\theta\|_4^2 \le a^3 \sqrt{n}\|\theta\|_\infty^2
  \le a^3 n \|\theta\|_\infty^3,$$
which proves the first statement. 
The randomized statement follows from
part \ref{I:Gordon-infty} of Lemma \ref{T:Gordon}, and the total variation
bound follows from Lemma \ref{T:tv}.
\end{enumerate}
\end{proof}

\bigskip

\begin{proof}[Proof of Corollary \ref{T:lp}]
The square negative correlation property was proved for the uniform measure 
on the ball of $\ell_p^n$ in \cite{BP} (see also \cite{ABP}) and for the
cone measure in \cite{NR}. The uniform measure on the ball is log-concave
by the Brunn-Minkowski theorem,
and it is not hard to show (cf.\ \cite{NR}) that $\E |X_i|^4 \le c$ for
some absolute constant in the case of cone measure.
Finally, it is well-known that if the
$\ell_p^n$ ball is scaled so that its uniform measure is isotropic, then it is
contained in $\big[-a_pn^{1/p},a_pn^{1/p}\big]^n$, where $a_p$ 
is a constant depending only on $p$;
it is not hard to show the same is true of the normalized cone measure.

Using all these facts, the statements for uniform measure on the ball
and cone measure on the sphere follow from Corollary \ref{T:SNCP};
the total variation bound for the ball follows from Lemma \ref{T:tv}.

The statements for the surface measure then follow from the fact,
proved in \cite{NR,Naor}
that the total variation distance between the cone and surface measures
is at most $\frac{c}{\sqrt{n}}$ for some absolute constant $c$, and both
of the error estimates are of at least this order (cf.\ the proofs of
\ref{I:log-concave} and \ref{I:bounded} of Corollary \ref{T:SNCP}).

\end{proof}

\bigskip

\begin{proof}[Proof of Corollary \ref{T:simplex}]
First, if $\Delta_n=\sqrt{n(n+2)}\conv\{v_1,\ldots,v_{n+1}\}$ is a 
regular simplex, then the $v_i$ form a normalized tight frame and also satisfy
\begin{equation}\label{E:0-sum}
\sum_{i=1}^{n+1}v_i=0
\end{equation}
Both of these facts are well-known and 
can be seen as consequences of John's theorem on
contact points between a convex body and the minimal volume ellipsoid 
containing it \cite{John}.   

To see Corollary \ref{T:simplex} as a consequence of Theorem \ref{T:proj}, 
consider the vectors
$$u_{ij}=\sqrt{\frac{n}{2(n+1)}}(v_i-v_j), \qquad 1\le i,j \le n+1, \quad
i\neq j.$$
It is not hard to show from \eqref{E:0-sum} that $\|u_{ij}\|_2=1$ for each
$i \neq j$, and that the $u_{ij}$ form a normalized tight frame because the $v_i$
do. Reflection in $u_{ij}^\perp$ is a reflection which interchanges
the vertices in the directions $v_i$ and $v_j$ and leaves all other vertices
of $\Delta_n$ fixed.
Theorem \ref{T:proj} can thus be applied (with $m=n(n+1)$), 
provided that $X$ is indeed isotropic under this scaling.  This is not 
obvious at this point, but follows easily from (\ref{E:moments}) below.

In order to compute the relevant expectations, one can embed $\Delta_n$ 
isometrically in $\R^{n+1}$ by the affine map with
$$\sqrt{n(n+2)}v_i \longmapsto \sqrt{(n+1)(n+2)}e_i;$$
the image of $\Delta_n$ under this map is 
$$\Delta_n'=\sqrt{(n+1)(n+2)}\conv\{e_1,\ldots,e_{n+1}\}.$$
Let $Y$ be the image of $X$ under this isometry; $Y$ is 
uniformly distributed in $\Delta_n'$.  Then
$$\inprod{X}{v_i}=\inprod{Y}{\sqrt{\frac{n+1}{n}}\left(e_i-\frac{1}{(n+1)}
\sum_{j=1}^{n+1}e_j \right)},$$
and so, adapting the notation of Theorem \ref{T:proj},
\begin{equation*}
\begin{split}
X_{(ij)} &= \inprod{X}{u_{ij}} \\
& = \sqrt{\frac{n}{2(n+1)}}\left(\inprod{X}{v_i}- \inprod{X}{v_j}\right)\\
&=\sqrt{\frac{n}{2(n+1)}}\left(\inprod{Y}{\sqrt{\frac{
n+1}{n}}e_i}-\inprod{Y}{\sqrt{\frac{n+1}{n}}e_j}\right)\\
&=\frac{1}{\sqrt{2}} \big(Y_i-Y_j\big).
\end{split}
\end{equation*}
The joint moments of the $Y_i$ are given by
\begin{equation}\label{E:moments}
\E\left[\prod_{i=1}^{n+1}Y_i^{r_i}\right]=\frac{\big[(n+1)(n+2)\big]^{r/2}
n!}{(n+r)!}\prod_{i=1}^{n+1}r_i!,
\end{equation}
where $r=\sum r_i.$ This formula follows easily from 
Lemma 1 of \cite{SZ}.
Using this one can show that $\E X_{(ij)}^2 = 1$ and thus that the
stated normalization for $\Delta_n$ is correct.
In addition, if $i\neq j$ and $k \neq l$, then
\begin{equation}\label{E:simplex-moments}
\E X_{(ij)}^2 X_{(kl)}^2 = \frac{(n+1)(n+2)}{(n+3)(n+4)} 
\begin{cases}
1 & \mbox{ if } \{i,j\} \cap \{k,l\} = \emptyset, \\
3 & \mbox{ if } |\{i,j\} \cap \{k,l\}| = 1, \\
6 & \mbox{ if } |\{i,j\} \cap \{k,l\}| = 2,
\end{cases}
\end{equation}
and 
\begin{equation} \label{E:simplex-3}
\E |X_{(ij)}|^3 \le 3\sqrt{2} \frac{\sqrt{(n+1)(n+2)}}{n+3} < 3\sqrt{2}.
\end{equation}
(What is really needed about this latter quantity is just that it is bounded 
by an absolute constant; this also follows immediately from Borell's lemma 
\eqref{E:Borell}.)

To estimate the first error term from Theorem \ref{T:proj}, by 
\eqref{E:simplex-moments},
\begin{equation*}
\begin{split}
\sum_{(ij),(kl)} & \theta_{(ij)}^2 \theta_{(kl)}^2 
  \E \bigl[X_{(ij)}^2 X_{(kl)}^2 \bigr] \\
&= \frac{(n+1)(n+2)}{(n+3)(n+4)}\left[ \sum_{(ij),(kl)} 
  \theta_{(ij)}^2 \theta_{(kl)}^2
  + 2 \sum_{|\{i,j\} \cap \{k,l\}|=1} \theta_{(ij)}^2 \theta_{(kl)}^2
  + 10 \sum_{(ij)} \theta_{(ij)}^4\right],
\end{split}
\end{equation*}
where in all of the above sums, indices run from $1$ to $n+1$, $i\neq j$,
and $k\neq l$. Since the $u_{ij}$ form a normalized tight frame,
$$\sum_{(ij)} \theta_{(ij)}^2 = \sum_{(ij)} \inprod{u_{(ij)}\otimes u_{(ij)}
  (\theta)}{\theta} = n+1,$$
and by the Cauchy-Schwarz inequality,
\begin{equation*}
\sum_{|\{i,j\} \cap \{k,l\}|=1} \theta_{(ij)}^2 \theta_{(kl)}^2
\le \sum_{|\{i,j\} \cap \{k,l\}|=1} \theta_{(ij)}^4
= 4 n \sum_{(ij)} \theta_{(ij)}^4.
\end{equation*}
By \eqref{E:0-sum} and the fact that the $v_i$ form a normalized tight frame,
\begin{equation*}
\begin{split}
\sum_{(ij)} \theta_{(ij)}^4 
&= \left(\frac{n}{2(n+1)}\right)^2 \sum_{i,j=1}^{n+1} 
  \big(\inprod{\theta}{v_i} -
  \inprod{\theta}{v_j}\big)^4 \\
& = \frac{n^2}{2(n+1)} \sum_{i=1}^{n+1} \inprod{\theta}{v_i}^4 + 6,
\end{split}
\end{equation*}
and so
\begin{equation}\label{E:simplex-1st}
\frac{1}{(n+1)^2} \sum_{(ij),(kl)} \theta_{(ij)}^2 \theta_{(kl)}^2 
  \E \bigl[X_{(ij)}^2 X_{(kl)}^2 \bigr] - 1
\le 4 \sum_{i=1}^{n+1} \inprod{\theta}{v_i}^4 + O\left(\frac{1}{n}\right).
\end{equation}

To estimate the second error term from Theorem \ref{T:proj},
\begin{equation}\label{E:simplex-2nd}
\begin{split}
\left(\sum_{(ij)} |\theta_{(ij)}|^3 \right)^{1/3} &=
  \sqrt{\frac{n}{2(n+1)}} \left(\sum_{i,j=1}^{n+1} 
  |\inprod{\theta}{v_i} - \inprod{\theta}{v_j}|^3 \right)^{1/3} \\
&\le \sqrt{\frac{n}{2}} \left(\sum_{i=1}^{n+1} |\inprod{\theta}{v_i}|^3
  \right)^3
\end{split}
\end{equation}
by the triangle inequality for the $\ell_3^{(n+1)^2}$ norm.
The first statement of the Corollary now follows by inserting 
\eqref{E:simplex-3}, \eqref{E:simplex-1st}, and \eqref{E:simplex-2nd} into
Theorem \ref{T:proj}, and noting as in the proof of part
\ref{I:log-concave} of Corollary \ref{T:SNCP}
that the first error term is of smaller order than the second error term.
The total variation bound then follows from Lemma \ref{T:tv}.

The randomized statement essentially follows from Lemma \ref{T:Gordon} as in
the proof of Corollary \ref{T:SNCP}. In this case one actually needs not
part \ref{I:Gordon-3} of Lemma \ref{T:Gordon}, but the same estimate for 
the norm
$$\| \theta \| = \left(\sum_{i=1}^{n+1} |\inprod{\theta}{v_i}|^3 
  \right)^{1/3},$$
which can be proved in the same way.
\end{proof}

\bigskip
{\bf Remarks:}
The same issue of estimating covariances of the squares
of frame components of $X$, which is explicit in the error bound of
Corollary \ref{T:uncon}, also arises implicitly in the proof of
Corollary \ref{T:simplex}.  In the latter case $\Cov(X_{(ij)}^2,X_{(kl)}^2)$
is not always small; the key point is that it is at worst a positive constant,
and this only happens for a negligible fraction of pairs $(ij),(kl)$.
Also, it is clear that the proof of Corollary \ref{T:simplex}
can be adapted to treat other distributions which posses the same 
symmetries as a centered regular simplex.

\section{Infinitesimal rotations} \label{S:rotations}

This section is mainly devoted to the proof of Theorem \ref{T:sph-symm}.
In the usual exchangeable pairs approach to Stein's method described in 
Section \ref{S:background}, one starts with a random variable $W$
and makes a small change to get a
new random variable $W'$ so that
the pair $(W,W')$ is exchangeable.  This pair is then used cleverly to 
estimate differences in expectations of test functions with respect to $W$ 
and some standard distribution.  Since the symmetries in Theorem 
\ref{T:sph-symm} are continuous rather than discrete, it is
possible to make an ``infinitesimal'' change in $W$ by making a 
small change, scaling appropriately, and then taking a limit.  

\begin{proof}[Proof of Theorem \ref{T:sph-symm}]
Suppose that $X$ is spherically symmetric and isotropic. By the spherical
symmetry, we may assume $\theta = e_1$, so $W=X_1$.

To prove the theorem, it suffices to bound $|\E f(W) - \E f(Z)|$, where 
$Z$ is a standard normal random variable and $f:\R \to \R$ is smooth with 
compact support. (The proof does not require this much 
regularity of $f$, but it produces no loss in generality.) 
Recall that the standard normal distribution is characterized by the 
identity
$$\E [g'(Z) - Z g(Z)] = 0$$
for all sufficiently regular $g$, and that given a test function $f$, 
$$g(t) = U_\nought f(t)=
e^{\frac{1}{2}t^2} \int_{-\infty}^t \bigl[f(s) - \E f(Z) \bigr] 
  e^{-\frac{1}{2}s^2}\ ds$$
satisfies the differential equation
$$g'(t) - t g(t) = f(t) - \E f(Z).$$
In particular
\begin{equation}\label{E:Stein-eq}
\E [g'(W) - W g(W)] = \E f(W) - \E f(Z).
\end{equation}

To carry out the infinitesimal rotations idea described above, define
a family of random variables $\{W_\epsilon\},$
 for $\epsilon \in (0,\frac{1}{2})$ as follows.
Let $A_\epsilon$ be the $n\times n$ orthogonal matrix
$$A_\epsilon = \begin{pmatrix} \sqrt{1-\epsilon^2} & \epsilon \\
  -\epsilon & \sqrt{1-\epsilon^2} \end{pmatrix}
  \oplus I_{n-2}.$$
Now let $U$ be a random $n\times n$ orthogonal matrix, chosen independently
of $X$ according to Haar measure; 
$U^TA_\epsilon U$ is a rotation in a random two-dimensional subspace through
an angle $\sin^{-1}(\epsilon)$.
Define
$$W_\epsilon = \inprod{(U^T A_\epsilon U) X}{e_1}.$$
By the rotational invariance of 
$X$, $(W,W_\epsilon)$ is an exchangeable pair for each $\epsilon$.

The following facts about the joint distribution of 
$(W,W_\epsilon)$ will be needed:
\begin{align}
\label{E:sph-symm-lemma-1}
\E[W-W_\epsilon|W] &= \bigl(1+O(\epsilon^2)\bigr)\dfrac{\epsilon^2}{n}W, \\
\label{E:sph-symm-lemma-2}
\E[(W-W_\epsilon)^2|W] 
  &= \bigl(1+O(\epsilon)\bigr) \dfrac{2\epsilon^2}{n} \E[X_2^2|W] + 
  O(\epsilon^4) W^2, \\
\label{E:sph-symm-lemma-3}
\E|W-W_\epsilon|^3 &= O(\epsilon^3).
\end{align}
Here and throughout this proof
the $O$ notation refers to asymptotic behavior as $\epsilon \to 0$, with
deterministic implied constants (that may depend on 
$n$, $f$, or the distribution of $X$). The proof of
\eqref{E:sph-symm-lemma-1} is given below. The proofs of
\eqref{E:sph-symm-lemma-2} and \eqref{E:sph-symm-lemma-3} are similar;
analogous estimates are proved in detail in \cite{EMeckes}.

First observe that by exchangeability,
$$\E \bigl((W- W_\epsilon) [g(W) + g(W_\epsilon)]\bigr) = 0,$$
because the expression is antisymmetric in $W$ and $W_\epsilon$.
(In the language of Section \ref{S:background}, this is essentially the 
observation that $\E T\alpha=0$, where $\alpha$ has been chosen as in the
remark at the end of Section \ref{S:background}.)  Now by Taylor's theorem,
\begin{equation*}
\begin{split}
(W - W_\epsilon) [g(W) + g(W_\epsilon)] &=
   2 (W - W_\epsilon) g(W) + (W - W_\epsilon)[g(W_\epsilon) - g(W)] \\
&=  2 (W - W_\epsilon) g(W) - (W - W_\epsilon)^2 g'(W) + R,
\end{split}
\end{equation*}
where
$$|R| \le \frac{1}{2}|W - W_\epsilon|^3 \|g''\|_\infty
  \le |W-W_\epsilon|^3 \| f' \|_\infty$$
by \eqref{E:U_o''}.

By \eqref{E:sph-symm-lemma-1}, \eqref{E:sph-symm-lemma-2},
\eqref{E:sph-symm-lemma-3}, and \eqref{E:U_o},
\begin{equation*}
\begin{split}
0 &= \frac{n}{2\epsilon^2} \E \bigl[2 g(W) (W-W_\epsilon) 
  - 2 (W - W_\epsilon)^2 g'(W) \bigr] + O(\epsilon) \\
&= \frac{n}{2\epsilon^2} \E \left[ \E \bigl[2 g(W) (W-W_\epsilon) 
  - 2 (W - W_\epsilon)^2 g'(W) \big| W \bigr] \right] + O(\epsilon) \\
&= \E \left[ W g(W) - \E[X_2^2 | W] g'(W) \right] + O(\epsilon),
\end{split}
\end{equation*}
and so, letting $\epsilon \to 0$, 
$$ \E [W g(W)] = \E \left[ \E[X_2^2 | W] g'(W) \right].$$
Therefore
\begin{equation}
\label{E:sph-symm-equality}
\E \bigl[f(W) - f(Z) \bigr] = \E [g'(W) - Wg(W)] 
= \E \left[ \bigl(1 - \E [X_2^2 | W]\bigr) g'(W)\right]
\end{equation}
for any smooth, bounded $f$. In particular,
$$d_{TV}\bigl(\Law(W), \gamma_\R \bigr) = \sup_f |\E f(W) - \E f(Z)|
  \le 4 \E \bigl| 1 - \E[X_2^2 | W] \bigr|,$$
where the supremum may be taken over smooth, compactly supported $f$
with $\| f \|_\infty \le 1$, so that $\| g' \|_\infty \le 4$ by \eqref{E:U_o'}.
This proves the first estimate.

To prove the second and third estimates, observe that by spherical symmetry,
$$\E \big[X_2^2 \big| X_1 \big]
  = \frac{1}{n-1} \left( \E \big[ \|X\|_2^2 \big| X_1\big] - X_1^2 \right),$$
and therefore
\begin{equation*}\begin{split}
\E \big| 1 - \E[X_2^2 | X_1 ] \big| &
  \le \frac{1}{n-1} \Big( \E \big| n - \E \big[\|X\|_2^2 \big| X_1 \big] \big|
   + \E \big| X_1^2 - 1\big| \Big) \\
&\le \frac{1}{n-1} \left(\E\big|\|X\|_2^2 -n \big|+2\right)\\
&\le \frac{1}{n-1} \left(\sqrt{\Var\left(\|X\|_2^2\right)}+2\right),
\end{split}\end{equation*}
by the Cauchy-Schwarz inequality and the isotropicity of $X$.

Finally, to prove \eqref{E:sph-symm-lemma-1}, first observe that
$$A_\epsilon = I_n + \left[\epsilon J - 
  \bigl(1 +O(\epsilon^2)\bigr) \frac{\epsilon^2}{2} I_2\right] 
\oplus 0_{n-2},$$
where $0_{n-2}$ denotes the $(n-2)\times (n-2)$ zero matrix and
$$J = \begin{pmatrix} 0 & 1 \\ -1 & 0 \end{pmatrix}.$$
Denote by $K$ the $2\times n$ matrix consisting of the first two rows of
the random orthogonal matrix $U$. Then
$$(U^T A_\epsilon U) X = X + K^T \left[\epsilon J - 
  \bigl(1 + O(\epsilon^2)\bigr) \frac{\epsilon^2}{2} I_2 \right]K X,$$
and so
\begin{equation}
\label{E:sph-symm-lemma-diff}
W-W_\epsilon = - \epsilon \inprod{(K^T J K)X}{e_1}
  + \bigl(1+O(\epsilon^2)\bigr) \frac{\epsilon^2}{2} \inprod{(K^T K)X}{e_1}.
\end{equation}
Now if $u_{ij}$ denote the entries of $U$, then by expanding in components,
$$
\E\left[\left.\inprod{(K^T J K)X}{e_1}\right|X\right] 
= \sum_{i=2}^n X_i \E (u_{11}u_{2i}-u_{21}u_{1i})
$$
and
$$
\E \left[\left.\inprod{(K^T K)X}{e_1}\right|X\right] 
= \sum_{i=1}^n X_i \E (u_{11}u_{1i} + u_{21}u_{2i}).
$$
Computing these expectations is not difficult because the distribution of
$U$ is unchanged by multiplying any row or column by $-1$, and any row or
column of $U$ is distributed uniformly on $S^{n-1}$. Therefore 
$$ \E u_{ij}u_{kl} = \delta_{ik}\delta_{jl} \frac{1}{n},$$
and so
$$\E\left[\left.\inprod{(K^T J K)X}{e_1}\right|X\right] = 0$$
and
$$\E \left[\left.\inprod{(K^T K)X}{e_1}\right|X\right] = \frac{2}{n}W.$$
Putting these together proves \eqref{E:sph-symm-lemma-1}.

The proofs of
\eqref{E:sph-symm-lemma-2} and \eqref{E:sph-symm-lemma-3} follow similarly
from \eqref{E:sph-symm-lemma-diff}. One needs in addition that
$\E |X_i|^3 < \infty$ and that $\E [X_i|X_1] = 0$ and $\E [X_i^2|X_1] = \E 
[X_2^2|X_1]$ for $i\neq 1$. One also needs values of fourth order moments
of the entries of $U$, which can be found, e.g., in \cite{AL}.

\end{proof}

\bigskip

The proof of Theorem \ref{T:unitary} is essentially the same as the 
proof of Theorem \ref{T:sph-symm}; the only difference is that $A_\epsilon$
is conjugated by a random unitary matrix instead of a random orthogonal 
matrix.  The relevant mixed moments of entries of a random unitary matrix 
can also be found in \cite{AL}.

\begin{proof}[Proof of Corollary \ref{T:Euclidean}]
If $X$ is uniformly distributed on the sphere of radius $\sqrt{n}$, then 
$$\Var \big( \|X\|_2^2 \big) = 0,$$
and the result follows immediately from
Theorem \ref{T:sph-symm}.
If $X$ is uniformly distributed on the ball of radius $\sqrt{n+2}$, then
it is easy to show by integration in polar coordinates that
$\Var \big( \|X\|_2^2 \big) < 4$, from which the result follows. In both
cases the constants, although not the order in $n$, can be improved by working
directly from the first estimate in Theorem \ref{T:sph-symm}.
\end{proof}

\bigskip

\begin{proof}[Proof of Corollary \ref{T:Poincare}]
By applying the Poincar\'e inequality \eqref{E:Poincare} to the function
$f(x) = \|x\|_2^2$, one obtains
$$\Var \big( \|X\|_2^2 \big) \le \frac{4n}{\lambda_1},$$
and so by Theorem \ref{T:sph-symm},
$$d_{TV} \big(\Law(W_\theta),\gamma_\R\big)
  \le \frac{8}{n-1} \left(\sqrt{\frac{n}{\lambda_1}} + 1\right).$$
By testing \eqref{E:Poincare} on a linear functional $f$, one obtains
that $\lambda_1 \le 1$ when $X$ is isotropic, and therefore
$$d_{TV} \big(\Law(W_\theta),\gamma_\R\big)
  \le \frac{8}{\sqrt{\lambda_1}(\sqrt{n}-1)}.$$
Since $d_{TV} \le 2$ always, this gives a result only for $n > 25$, for which
the stated estimate now follows.
\end{proof}

\bibliographystyle{plain}
\bibliography{clt}

\begin{thebibliography}{10}

\bibitem{ABP}
M.~Anttila, K.~Ball, and I.~Perissinaki.
\newblock The central limit problem for convex bodies.
\newblock {\em Trans. Amer. Math. Soc.}, 355(12):4723--4735 (electronic), 2003.

\bibitem{AL}
S.~Aubert and C.~S. Lam.
\newblock Invariant integration over the unitary group.
\newblock {\em J. Math. Phys.}, 44(12):6112--6131, 2003.

\bibitem{BP}
K.~Ball and I.~Perissinaki.
\newblock The subindependence of coordinate slabs in {$l\sp n\sb p$} balls.
\newblock {\em Israel J. Math.}, 107:289--299, 1998.

\bibitem{BB}
J.~Bastero and J.~Bernu{\'e}s.
\newblock Asymptotic behaviour of averages of $k$-dimensional margins of
  measures on {$\R^n$}.
\newblock Preprint, 2005.

\bibitem{Bobkov3}
S.~G. Bobkov.
\newblock Personal communication.

\bibitem{Bobkov1}
S.~G. Bobkov.
\newblock On concentration of distributions of random weighted sums.
\newblock {\em Ann. Probab.}, 31(1):195--215, 2003.

\bibitem{Bobkov4}
S.~G. Bobkov.
\newblock Spectral gap and concentration for some spherically symmetric
  probability measures.
\newblock In {\em Geometric Aspects of Functional Analysis (2001--2002)},
  volume 1807 of {\em Lecture Notes in Math.}, pages 37--43. Springer, Berlin,
  2003.

\bibitem{BK}
S.~G. Bobkov and A.~Koldobsky.
\newblock On the central limit property of convex bodies.
\newblock In {\em Geometric Aspects of Functional Analysis (2001--2002)},
  volume 1807 of {\em Lecture Notes in Math.}, pages 44--52. Springer, Berlin,
  2003.

\bibitem{BHVV}
U.~Brehm, P.~Hinow, H.~Vogt, and J.~Voigt.
\newblock Moment inequalities and central limit properties of isotropic convex
  bodies.
\newblock {\em Math. Z.}, 240(1):37--51, 2002.

\bibitem{BV}
U.~Brehm and J.~Voigt.
\newblock Asymptotics of cross sections for convex bodies.
\newblock {\em Beitr\"age Algebra Geom.}, 41(2):437--454, 2000.

\bibitem{Bryc}
W.~Bryc.
\newblock {\em The Normal Distribution, Characterizations with Applications},
  volume 100 of {\em Lecture Notes in Statistics}.
\newblock Springer-Verlag, New York, 1995.

\bibitem{CDM}
S.~Chatterjee, P.~Diaconis, and E.~Meckes.
\newblock Exchangeable pairs and {P}oisson approximation.
\newblock {\em Probab. Surv.}, 2:64--106 (electronic), 2005.

\bibitem{DF2}
P.~Diaconis and D.~Freedman.
\newblock A dozen de {F}inetti-style results in search of a theory.
\newblock {\em Ann. Inst. H. Poincar\'e Probab. Statist.}, 23(2,
  suppl.):397--423, 1987.

\bibitem{Gordon}
Y.~Gordon.
\newblock On {M}ilman's inequality and random subspaces which escape through a
  mesh in {${\bf R}\sp n$}.
\newblock In {\em Geometric aspects of functional analysis (1986/87)}, volume
  1317 of {\em Lecture Notes in Math.}, pages 84--106. Springer, Berlin, 1988.

\bibitem{HR}
S.~Holmes and G.~Reinert.
\newblock Stein's method for the bootstrap.
\newblock In {\em Stein's Method: Expository Lectures and Applications},
  volume~46 of {\em IMS Lecture Notes Monogr. Ser.}, pages 95--136. Inst. Math.
  Statist., Beachwood, OH, 2004.

\bibitem{John}
F.~John.
\newblock Extremum problems with inequalities as subsidiary conditions.
\newblock In {\em Studies and Essays Presented to R. Courant on his 60th
  Birthday, January 8, 1948}, pages 187--204. Interscience Publishers, Inc.,
  New York, N. Y., 1948.

\bibitem{KLS}
R.~Kannan, L.~Lov{\'a}sz, and M.~Simonovits.
\newblock Isoperimetric problems for convex bodies and a localization lemma.
\newblock {\em Discrete Comput. Geom.}, 13(3-4):541--559, 1995.

\bibitem{KL}
A.~Koldobsky and M.~Lifshits.
\newblock Average volume of sections of star bodies.
\newblock In {\em Geometric Aspects of Functional Analysis (1996--2000)},
  volume 1745 of {\em Lecture Notes in Math.}, pages 119--146. Springer,
  Berlin, 2000.

\bibitem{Ledoux2}
M.~Ledoux.
\newblock Spectral gap, logarithmic {S}obolev constant, and geometric bounds.
\newblock {\em J. Differential Geom.}
\newblock To appear.

\bibitem{Ledoux1}
M.~Ledoux.
\newblock {\em The Concentration of Measure Phenomenon}, volume~89 of {\em
  Mathematical Surveys and Monographs}.
\newblock American Mathematical Society, Providence, RI, 2001.

\bibitem{EMeckes}
E.~Meckes.
\newblock Manuscript in preparation.

\bibitem{Naor}
A.~Naor.
\newblock The surface measure and cone measure on the sphere of {$\ell_p^n$}.
\newblock {\em Trans. Amer. Math. Soc.}
\newblock To appear.

\bibitem{NR}
A.~Naor and D.~Romik.
\newblock Projecting the surface measure of the sphere of {$\ell_p^n$}.
\newblock {\em Ann. Inst. H. Poincar\'e Probab. Statist.}, 39(2):241--261,
  2003.

\bibitem{Prekopa}
A.~Pr{\'e}kopa.
\newblock On logarithmic concave measures and functions.
\newblock {\em Acta Sci. Math. (Szeged)}, 34:335--343, 1973.

\bibitem{Reitzner}
M.~Reitzner.
\newblock Central limit theorems for random polytopes.
\newblock {\em Probab. Theory Related Fields}.
\newblock To appear.

\bibitem{RR}
Y.~Rinott and V.~Rotar.
\newblock On coupling constructions and rates in the {CLT} for dependent
  summands with applications to the antivoter model and weighted
  {$U$}-statistics.
\newblock {\em Ann. Appl. Probab.}, 7(4):1080--1105, 1997.

\bibitem{SZ}
G.~Schechtman and J.~Zinn.
\newblock On the volume of the intersection of two {$L\sp n\sb p$} balls.
\newblock {\em Proc. Amer. Math. Soc.}, 110(1):217--224, 1990.

\bibitem{Sodin}
S.~Sodin.
\newblock Tail-sensitive {G}aussian asymptotics for marginals of concentrated
  measures in high dimension.
\newblock Preprint, 2005.

\bibitem{Stein1}
C.~Stein.
\newblock {\em Approximate {C}omputation of {E}xpectations}, volume~7 of {\em
  Institute of Mathematical Statistics Lecture Notes---Monograph Series}.
\newblock Institute of Mathematical Statistics, Hayward, CA, 1986.

\bibitem{Stein2}
C.~Stein.
\newblock The accuracy of the normal approximation to the distribution of the
  traces of powers of random orthogonal matrices.
\newblock Technical Report 470, Stanford University Dept. of Statistics, 1995.

\bibitem{Sudakov}
V.~N. Sudakov.
\newblock Typical distributions of linear functionals in finite-dimensional
  spaces of high dimension.
\newblock {\em Dokl. Akad. Nauk SSSR}, 243(6):1402--1405, 1978.

\end{thebibliography}

\end{document}